\documentclass[runningheads]{llncs}

\usepackage{algorithm,algpseudocode,amsfonts,amsmath,amssymb} 
\usepackage{adjustbox,booktabs,graphicx,hyperref,makecell,multirow,siunitx,wrapfig,xparse,xspace} 
\usepackage{caption} 
\usepackage{subcaption}
\usepackage{arydshln} 
\usepackage[shortlabels,inline]{enumitem} 
\usepackage[numbers]{natbib} 
\PassOptionsToPackage{table,xcdraw}{xcolor}

\graphicspath{{figures/pdc/}{figures/polyhedra/}}

\usepackage{definitions/defs}  
\usepackage[nameinlink]{cleveref} 

\author{Shannon Kelley\inst{1} \and Aleksandr M. Kazachkov\inst{2} \and Ted Ralphs\inst{1}}
\institute{Lehigh University, Bethlehem, PA, USA \and University of Florida, Gainesville, FL, USA}
\date{\today}

\title{Parametric Disjunctive Cuts for Sequences of Mixed Integer Linear Optimization Problems}
\titlerunning{Parametric Disjunctive Cuts}

\begin{document}
    \maketitle

    \begin{abstract}

        Many applications require solving sequences of related mixed-integer linear programs. We introduce a class of \emph{parametric disjunctive inequalities} (PDIs), obtained by reusing the disjunctive proofs of optimality from prior solves to construct cuts valid for  perturbed instances. We describe several methods of generating such cuts that navigate the tradeoff between computational expense and strength. We provide sufficient conditions under which PDIs support the disjunctive hull and a tightening step that guarantees support when needed. On perturbed instances from MIPLIB 2017, augmenting branch-and-cut with PDIs substantially improves performance, reducing total solve times on the majority of challenging cases.

    \end{abstract}
	\keywords{parametric optimization, warm starting, disjunctive cuts, MILP}

    \section{Introduction}
    \label{sec:intro}


	We propose a class of \emph{parametric disjunctive inequalities} (PDIs) for accelerating the solution of sequences of mixed integer linear optimization problems (MILPs). There are many methodologies whose implementation requires solving sequences of instances from a parametric family of MILPs, including column generation, Lagrangian relaxation, Benders decompositions, mixed integer nonlinear optimization, bilevel optimization, multiobjective optimization, and online optimization \cite{Barnhart98branch-and-price,Beasley93lagrangian-relaxation,Hassanzadeh14generalized-benders,Duran86outer-approximation,Tahernejad20bilevel,Fallah2024value-function,Nair20neural-diving}. Applications are diverse, including energy problems, supply chain optimization, and revenue management \cite{Xavier21learning-unit-commitment,Nazari18rl-for-vrp,Gasse19learn2branch}.


	The similar structure of the instances in parametric families makes them natural candidates for warm starting. Existing approaches reuse solutions and disjunctions, transfer parameters and branching policies, guide cut selection \cite{Nair20neural-diving,Guzelsoy09dual-methods,Munoz23compressing-mip-trees,Duran86outer-approximation,Gasse19learn2branch,Deza23ml-cut-survey}, or even re-use the entire branch-and-bound tree~\cite{Guzelsoy09dual-methods}. Yet, despite the central role of cutting planes in performance \cite{Bixby07mip-progress}, the cut generation process itself is rarely warm started.

    
    This omission is understandable: most cuts are inexpensive to regenerate, so parameterization can appear unnecessary. Two observations motivate a different view. First, cuts are typically generated with default settings, and recent work shows that learning instance-dependent parameters can improve solver performance~\cite{Becu24approximating-gomory,Chetelat23continuous-cut-algorithm,Clarke23differentiable-cutting-plane}. Second, some families---especially disjunctive cuts from large disjunctions---offer among the strongest relaxations~\cite{Perregaard01prlp,Chen11disjunctive-programming-characterizations,Kazachkov22vpc}, yet are often disabled in practice due to their computational cost~\cite{Bixby07mip-progress,Cbc}. We address how to obtain strong disjunctive cuts without repeatedly solving the auxiliary LPs they typically require: our PDIs exploit shared structure to amortize generation across related solves, reducing time while retaining much of the original strength.

    


	The validity of the disjunctions we construct is based solely on integrality of the variables, arising as the leaf nodes of (partial) \textit{branch-and-bound} trees constructed using the standard strategy of branching on variables. The key insight behind parametric disjunctive inequalities is that the certificate of validity of an inequality generated from such a disjunction can be easily exploited to provide an inequality and a certificate of validity for \emph{any} other MILP with the same number of constraints and the same integer variables. We focus on cuts derived from disjunctions with a large number of terms, which are expensive to generate but powerful. We show that parameterizing their certificates across a sequence of instances effectively amortizes generation cost. We then evaluate whether this accelerates overall solution time of structurally similar problems.


	\paragraph{Contributions.}
	\begin{enumerate}[i)]
		\item We introduce a class of PDIs for which we describe how to parameterize their certificates of validity and use these parametric certificates to generate inequalities valid for new instances.
		\item We prove conditions under which a PDI, when applied to a new instance, supports the disjunctive hull, as well as introduce a tightening procedure that ensures this property holds.
		\item We integrate the generation of PDIs into branch and cut and evaluate performance on families of instances derived by perturbing MIPLIB 2017 instances. We find consistent improvements over the generation procedure described in~\citet{Kazachkov22vpc} (which is not parametric) and over a default baseline that does not use disjunctive cuts for challenging cases.
	\end{enumerate}

    \section{Notation and Preliminaries}
    \label{sec:background_ipco}

	Our goal is to generate inequalities valid for an MILP drawn from a parametric family $\MIPSet$ of related instances. While our techniques apply as long as the set of integer-constrained variables is consistent across instances, we assume for simplicity that all instances are defined over the same variables and number of constraints.

	Each instance $\MIPk \in \MIPSet$ takes the form:
	\begin{equation}
		\min_{x \in \PIk} \quad c^k x
		\tag{IP-$\MIPk$} \label{IP-k}
	\end{equation}
	where $\PIk \defeq \left\{x \in \Pk : x_j \in \Z \text{ for all } j \in \intvars \right\} \text{ and } \Pk \defeq \left\{x \in \R^n : A^k x \ge b^k \right\},$
	where $A^k \in \Q^{\numRowsP \times n}$, $b^k \in \Q^{\numRowsP}$, $c^k \in \Q^{1 \times n}$, and $\intvars \subseteq [n] \defeq \{1,\ldots,n\}$ is the set of indices of integer variables. We assume the constraints include nonnegative upper and lower bounds on all variables. 
    For matrix $A$, we represent its $i^\text{th}$ row as $\mxrow{A}{i}$ and its $j^\text{th}$ column as $\mxcol{A}{j}$.

	To strengthen the LP relaxation, we rely on \emph{disjunctive cuts}, which are valid inequalities derived from valid disjunctions.

	\begin{definition}
		A \emph{valid inequality} for a set $\S \subseteq \Re^n$ is a pair $(\alpha, \beta) \in \Re^{n} \times \R$ such that $\S \subseteq \left\{x \in \Re^n \suchthat \alpha^\T x \ge \beta \right\}.$
	\end{definition}
	\begin{definition}
		A \emph{valid disjunction} for a set $\S \subseteq \Re^n$ is a collection $\disj \defeq \disjSet$, where $\disjTermsIndexSet$ is an index set, $\Xt \subset \Re^n$ for $t \in \disjTermsIndexSet$, and $\S \subseteq \X$.
	\end{definition}
	In this work, we assume $\Xt \defeq \left\{x \in \R^n \suchthat D^t x \ge D^t_0 \right\},$ where $D^t \in \Q^{\numRowsDkt \times n}$ and $D^t_0 \in \Q^{\numRowsDkt}$. In fact, we further assume $\Xt$ is a hyperrectangle obtained by restricting the bounds on one or more variables. 

	We frequently reference the feasible region of an instance $k \in \MIPSet$ restricted to a disjunctive term $t \in \disjTermsIndexSet$:
	\begin{equation*} 	A^{\MIPk t} \defeq
	\begin{bmatrix}
		A^k \\
		D^t
	\end{bmatrix}, \quad
	b^{\MIPk t} \defeq
	\begin{bmatrix}
		b^k \\
		D^t_0
	\end{bmatrix}, \quad \text{and} \quad
	\Qkt \defeq \Pk \cap \Xt = \left\{x \in \R^n : A^{\MIPk t} x \ge b^{\MIPk t} \right\}. 
    \end{equation*}
	Additionally, we refer to the convex hull of feasible regions defined by all terms of a disjunction as the \emph{disjunctive hull}, defined as $ \PDk \defeq \cl\conv\left(\bigcup_{t \in \disjTermsIndexSet} \Qkt \right)$, and we index the subset of disjunctive terms that are feasible for \IP{k} as 
    $
        \feasibleDisjTermsIndexSetk
        \defeq \{t \in \disjTermsIndexSet : \Qkt \neq \emptyset\}.
    $

	We refer to the set of indices of any collection of $n$ linearly independent constraints binding at some basic solution as 
    a \emph{basis}. Conceptually, this is equivalent to the concept of a basis that arises in the theory of LPs, when the problem is usually assumed to be in standard form. Bases in standard form are in one-to-one correspondence to what we are here calling a basis and as such, they can be feasible or infeasible, as specified in the following definition.
%
	\begin{definition}
		A basis $\activeConstraintIndexSett$ is \emph{feasible} for $\Qkt$ if the associated basic solution $(\mxrow{A^{kt}}{\activeConstraintIndexSett})^{-1} b^{kt}_{\activeConstraintIndexSett}$ lies in $\Qkt$.
	\end{definition}
%


    \paragraph{Disjunctive Cut Generation.}

	Valid inequalities for the disjunctive hull $\PDk$ can be generated by solving an LP, either the classical \emph{Cut Generating LP (CGLP)}~\cite{Balas79disjunctive-programming,Chen11disjunctive-programming-characterizations} or the so-called \emph{Point-Ray LP (PRLP)}~\cite{Perregaard01prlp,Kazachkov22vpc}. Regardless of which of these two LPs is used, the solution is an inequality $\cut \in \Re^n \times \R$ that satisfies:
	\[
		\tag{FL} \label{farkas}
		\left.
        \begin{array}{l}
        	{\alpha} ^\T = \vt A^{kt}     \\
            \beta \le \vt b^{kt} \\
            \vt \in \nonnegreals^{1 \times (\numRowsQkt)}
        \end{array}
        \ \ \right\}
        \ \
        \rlaptext{for all $t \in \disjTermsIndexSet$.}
	\]
	By \emph{Farkas' Lemma}~\cite{Farkas02}, this condition is equivalent to the inequality $\cut$ being valid for $\PDk$. We refer to $\singleCertificate$ as the \emph{Farkas multipliers} or \emph{Farkas certificate} for $\cut$ relative to the disjunction $\disj$ and instance \IP{k}.

	The method used to compute $\singleCertificate$ depends on whether we are solving a CGLP or a PRLP. In the CGLP, condition~\eqref{farkas} is included directly in the formulation, making $\singleCertificate$ part of the solution. In contrast, the PRLP reduces the dimensionality of the problem by eliminating the variables whose values produce the Farkas certificate~\cite{Perregaard01prlp}, allowing for more efficient computation with larger disjunctions and, consequently, stronger cuts~\cite{Kazachkov22vpc}. For this reason, we adopt the PRLP when generating disjunctive cuts from an LP. The corresponding Farkas certificate is then reconstructed post hoc, as described by \citet[Lemma 3]{Kazachkov23strengthen}. In practice, each term $\vt$ is computed from $\activeConstraintIndexSett$ and $\Qkt$, selected to represent the closest feasible basis to the cut $\cut$ for term $t$. If $\Qkt$ is infeasible and no such basis exists, we assign $\vt = 0$ by default.

    \section{Parametric Disjunctive Cuts}
   \label{sec:PDI}

	\paragraph{Certificates.} The certificate of validity for a PDI consists of a disjunction and a set of multipliers (the Farkas certificate) associated with each term of the disjunctions that together satisfy~\eqref{farkas}, certifying validity of the cut.

	\begin{definition}
    		Let $\MIPk \in \MIPSet$, $\singleCertificate$, and a disjunction $\{\Xt\}_{t \in \disjTermsIndexSet}$ valid for \IP{k} be given. Then $(\singleCertificate, \{\Xt\}_{t \in \disjTermsIndexSet}$) is a \emph{certificate of validity} for $(\alpha^\T, \beta)$ if 
		$
		\alpha_j \defeq
		\max_{t \in \disjTermsIndexSet}
		\{ v^t \mxcol{A^{\ell t}}{j} \}
		$
		and
		$
		\beta \defeq
		\min_{t \in \disjTermsIndexSet}
		\{ v^t b^{\ell t} \}
		$
		for all $ j \in [n] $.
	\end{definition}
%
%
    In our algorithms, certificates of validity are constructed with respect to one particular instance in a family and then are used to generate cuts for other instances. 
	\begin{definition}
		An MILP instance \IP{k} and disjunction $\disjSet$ \emph{induce} the Farkas certificate $\vt$ if $\vt[t_1] A^{kt_1} = \vt[t_2] A^{kt_2}$ for all $(t_1, t_2) \in \feasibleDisjSetk \times \feasibleDisjSetk$.
	\end{definition}
	To identify the basis used to construct the Farkas coefficients for a given feasible disjunctive term $\Qkt$, we associate $\vt$ with $\activeConstraintIndexSett$, a basis of $\Qkt$, which contains the indices of the constraints whose linear combination produces the cut coefficients~$\alpha$.

	\begin{definition}
		The basis $\activeConstraintIndexSett$ \emph{determines} the Farkas certificate $\vt$ if $\{i \in [\numRowsQkt] : \vt_i > 0\} \subseteq \activeConstraintIndexSett$.
	\end{definition}

    \paragraph{Parametric Cut Generation.} To reuse disjunctive inequalities across related instances while preserving validity, we derive parametric inequalities (parametric cuts) for the disjunctive relaxation via a procedure that guarantees their validity by ensuring such validity is satisfied by a known certificate. 

    \begin{definition}
		Let $\MIPSet$ be the indices of a family of MILPs, $\{\Xt\}_{t \in \disjTermsIndexSet}$ be a disjunction valid for every member of the family, and  $\singleCertificate$ be a Farkas certificate satisfying~\eqref{farkas}. Then $\{(\alpha^k, \beta^k)\}_{k \in \MIPSet}$ is a \emph{family of parametric disjunctive inequalities (PDIs)} with respect to $(\singleCertificate, \{\Xt\}_{t \in \disjTermsIndexSet})$ if for all $k \in \MIPSet$, we have that $(\singleCertificate, \{\Xt\}_{t \in \disjTermsIndexSet})$ is a certificate of validity of $(\alpha^k, \beta^k)$ for $\PIk$. 
	\end{definition}

    

    \paragraph{Validity.}
    \label{subsec:validity}

\begin{figure}[htbp]
		\centering
		\begin{subfigure}{0.475\textwidth}
			\centering
			\includegraphics[width=\textwidth]{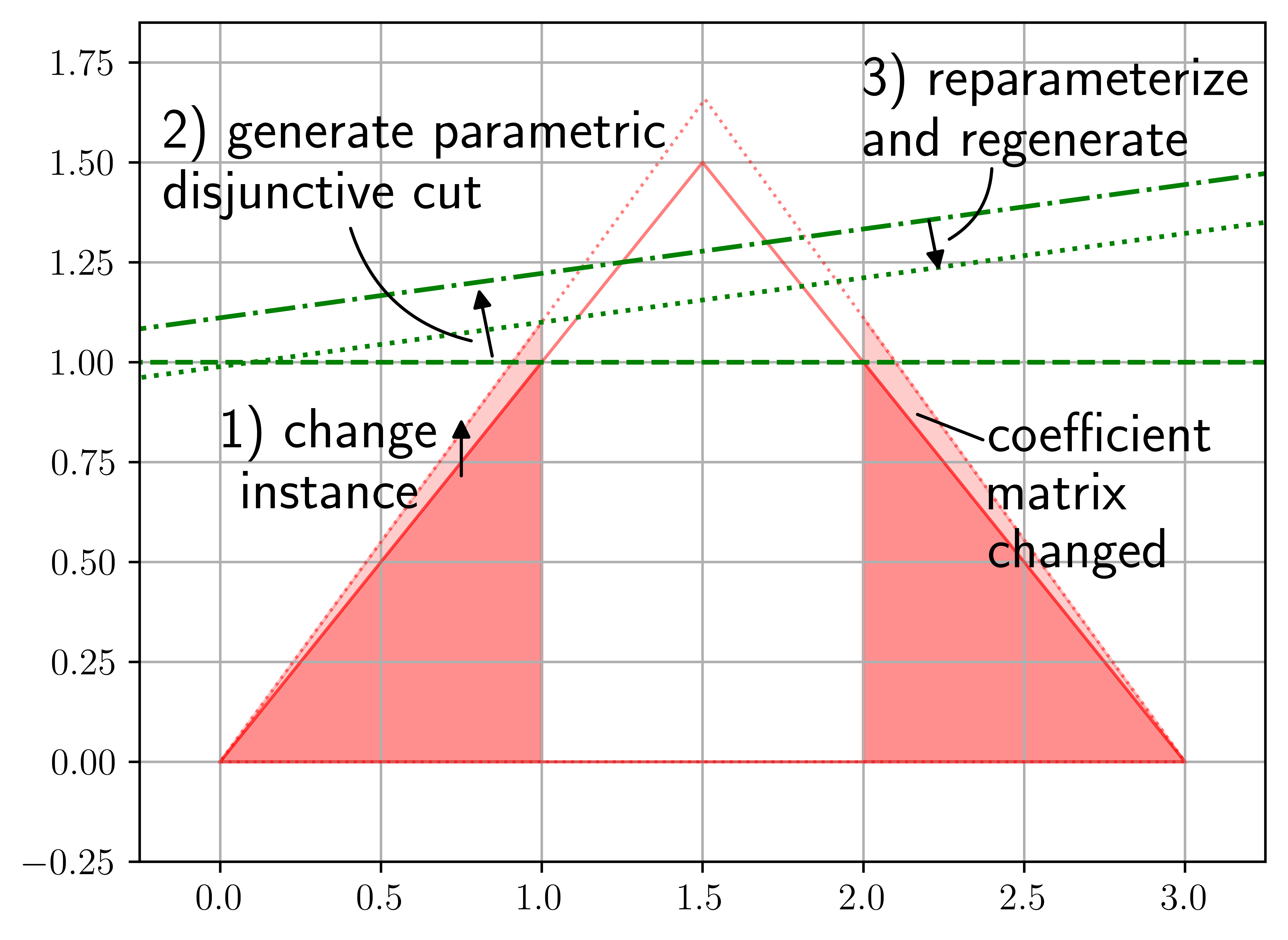}
			\caption{matrix perturbation}
			\label{p:tight_matrix}
		\end{subfigure}
		\begin{subfigure}{0.475\textwidth}
			\centering
			\includegraphics[width=\textwidth]{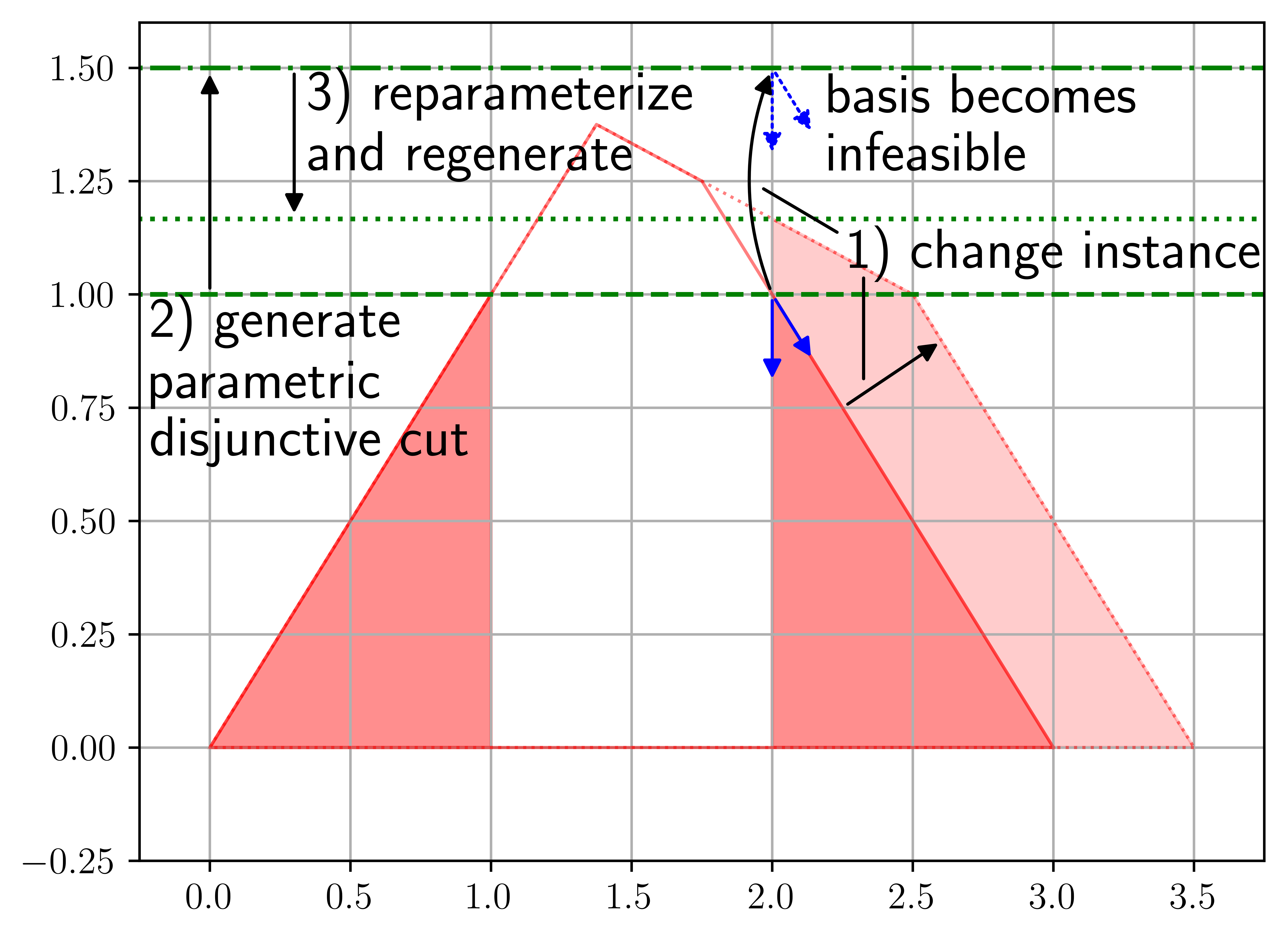}
			\caption{infeasibility from perturbation}
			\label{p:tight_basis}
		\end{subfigure}
		\caption{Examples showing that PDIs from the same certificate need not all support the disjunctive hull}
	\end{figure}

	We now establish procedures for generating parametric families of valid inequalities (proofs of correctness in \Cref{a:proofs}) by taking a previously generated certificate of validity as input and generating a cut satisfying~\eqref{farkas} from it. By taking the Farkas certificate as an input, we avoid solving solving an LP to compute it, while hopefully retaining the strength of the cuts due to the similarity between the instances.
    \begin{definition}
        \label{def:farkas_family}
		A \emph{family of Farkas PDIs} is $\{(\alpha^k, \beta^k)\}_{k \in \MIPSet}$ such that $(\alpha^k, \beta^k) \defeq ([\max_{t \in \disjTermsIndexSet} \{ v^t \mxcol{A^{k t}}{1}, ...,  \max_{t \in \disjTermsIndexSet} \{ v^t \mxcol{A^{kt}}{n}\}]^{\T}, \min_{t \in \disjTermsIndexSet}
		\{ v^t b^{kt} \})$.
	\end{definition}  
    For ease of exposition, when we describe families of PDIs going forward, we mean \emph{Farkas} PDIs. The following result shows they are parametrically valid for $\PDk[\ell]$. For clarity, we use the indices $\MIPk$ and $\ell$ to denote, respectively, the instance that induces the Farkas certificate and the instance to which it is applied.
	\begin{theorem}
		\label{t:validate-parameterization}
		Let $\MIPk, \ell \in \MIPSet$. Then given a Farkas certificate $\singleCertificate$ induced by disjunction $\{\Xt\}_{t \in \disjTermsIndexSet}$ and \IP{k}, we have that $\alpha^\T x \ge \beta$ for all $x \in \PDk[\ell]$, where
		$
		\alpha_j \defeq
		\max_{t \in \disjTermsIndexSet}
		\{ v^t \mxcol{A^{\ell t}}{j} \}
		$
		and
		$
		\beta \defeq
		\min_{t \in \disjTermsIndexSet}
		\{ v^t b^{\ell t} \}
		$
		for all $ j \in [n] $.
	\end{theorem}
    \Cref{p:tight_matrix,p:tight_basis} illustrate cuts produced by \Cref{t:validate-parameterization}, marked in the figures with the annotations ``2) generate parametric disjunctive cut''.
    \paragraph{Strengthening.}
%
%
	A natural question arises when generating a PDI $\cut$ for a disjunction $\disjSet$ and instance \IP{k}: does the cut support the disjunctive hull $\PDk$?
	\begin{definition}
		The inequality $\cut$ valid for a polyhedron $\mathcal{P}$ \emph{supports} $\mathcal{P}$ if there exists $x \in \mathcal{P}$ such that $\alpha^{\T} x = \beta$.
	\end{definition}
	Cuts generated by solving the system~\eqref{farkas} are guaranteed to support the disjunctive hull, but PDIs from \Cref{t:validate-parameterization} need not. The following conditions guarantee support, however.
	\begin{lemma}
        \label{lem:tighten_parameterization}
        Let $k, \ell \in \MIPSet$; 
        $\singleCertificate$ be a Farkas certificate induced by \IP{k} and disjunction $\disjSet$;
        $\{(\alpha^k, \beta^k)\}_{k \in \MIPSet}$ be the associated family of Farkas PDIs;
        and
        $\activeConstraintIndexSett$ be the basis determining $\vt$ for all $t \in \disjTermsIndexSet$.
        Then $(\alpha^\ell, \beta^\ell)$ supports $\PDk[\ell]$ if 
        \begin{enumerate}
            \item $\singleCertificate$ is induced by \IP{\ell} and $\disjSet$, and
            \item $\activeConstraintIndexSett$ is feasible for $\Qkt[\ell]$ for all $t \in \disjTermsIndexSet$.
        \end{enumerate}
    \end{lemma}
	Condition 1 fails when the coefficient matrix changes or a term that is infeasible in \IP{k} becomes feasible in \IP{\ell}. Condition 2 breaks when $\activeConstraintIndexSett$ becomes infeasible for some $\Qkt[\ell]$ when it was feasible for \IP{k}. The latter is often mitigated by removing infeasible terms before parameterization, leaving only bases that turn infeasible for still-feasible terms. See ``2)'' in \Cref{p:tight_matrix,p:tight_basis} below and the emergence of a newly feasible term in \Cref{p:tight_term} in \Cref{a:tables_figures}.
%
%

	The conditions of \Cref{lem:tighten_parameterization} hint at a method of recovering disjunctive hull support when it is lost. We recover support in two steps: (1) \emph{reparameterize} by computing the certificate for each term as needed, then (2) \emph{regenerate} the PDI using the updated certificates. The result of these steps are annotated ``3)'' in each of the three failure cases illustrated \Cref{p:tight_matrix,p:tight_term,p:tight_basis}.

	\paragraph{Step 1: Reparameterize to Find Supporting Certificates.}
	We refer to the process of obtaining a supporting cut and its corresponding \emph{supporting certificate}---the Farkas certificate that generates the supporting cut---as \emph{reparameterization}. \Cref{lem:supporting_certificate_exists} shows that this involves solving the primal-dual pair \ref{LP-lt-alpha} and \ref{Dual-lt-alpha}
	\begin{center}
		\begin{minipage}[t]{0.45\textwidth}
			\begin{equation}
				\begin{aligned}
					\min_{x \in \Qkt[\ell]} \; \alpha^{\T} x
				\end{aligned}
				\label{LP-lt-alpha} \tag{LP-$\ell t(\alpha)$}
			\end{equation}
		\end{minipage}
		\hfill
		\begin{minipage}[t]{0.45\textwidth}
			\begin{equation}
				\begin{aligned}
					\max_{v \in \Dlta} \; v b^{\ell t}
				\end{aligned}
				\label{Dual-lt-alpha} \tag{Dual-$\ell t(\alpha)$}
			\end{equation}
		\end{minipage}
	\end{center}
	where $\Dlta \defeq \{v \in \R^{1 \times (\numRowsQkt)} : v A^{\ell t} = \alpha^\T, \; v \geq 0\}$.
	\begin{lemma}
		\label{lem:supporting_certificate_exists}
		Let $\ell \in \MIPSet$, $ t \in \disjTermsIndexSet $, and $ \alpha \in \R^n $. If $\Qkt[\ell]$ is bounded and nonempty, then $(\alpha, \alpha^{\T} \bar{x})$ such that $\opt{x} \in \arg\min_{x \in \Qkt[\ell]}\{\alpha^{\T} x\}$ is a supporting inequality for $\Qkt$ with supporting certificate $\barvt \in \arg\max_{\vt \in \Dlta}\{\vt b^{\ell t}\}$.
	\end{lemma}
	Further, provided there is no constraint in $\Qkt[\ell]$ parallel to $\alpha$, \Cref{cor:supporting_certificate_unique} shows that the only supporting cut with coefficients $\alpha$ and a nonnegative Farkas certificate is the one derived from the optimal solutions to the primal–dual pair \ref{LP-lt-alpha} and \ref{Dual-lt-alpha}. This is essential, since \Cref{t:validate-parameterization} guarantees the validity of PDIs only when the Farkas certificate is nonnegative.
	\begin{corollary}
		\label{cor:supporting_certificate_unique}
		Let $\ell \in \MIPSet$, $ t \in \disjTermsIndexSet $, and $ \alpha \in \R^n $. If $\Qkt[\ell]$ is bounded and nonempty, and $ \alpha \neq \mxrow{A^{\ell t}}{i} $ for all $ i \in [\numRowsQkt] $, then there exists unique $\opt{x} \in \Qkt[\ell]$ and unique $\barvt \in \Dlta$ such that $ \alpha^{\T} \opt{x} = \barvt b^{\MIPk t} $.
	\end{corollary}
	Because the supporting certificate is typically unique for fixed $(\alpha,t)$, solving \ref{LP-lt-alpha} is unavoidable; however, storing the basis that originally determined each certificate often accelerates reoptimization after small perturbations. The procedure is agnostic to which condition in \Cref{lem:tighten_parameterization} failed---it simply recovers, per term, a certificate that supports the given $\alpha$. This is precisely the result that, as we show next, is required for \Cref{t:validate-parameterization} to produce PDIs that support the disjunctive hull.
	\paragraph{Step 2: Regenerate to Yield Supporting PDIs.}
	We leverage \Cref{lem:supporting_certificate_exists} to construct a Farkas certificate that satisfies the conditions of \Cref{lem:tighten_parameterization}, ensuring that \Cref{t:validate-parameterization} produces a cut supporting the disjunctive hull. This process is detailed in \Cref{alg:strong_param_vpcs}.
	\begin{algorithm}[t]
		\caption{Strong PDI Generator}
		\label{alg:strong_param_vpcs}
		\begin{algorithmic}[1]
			\Require $\MIPk, \ell, \disjSet, \singleCertificate, \activeConstraintIndexCollection$
			\Ensure \IP{k} and $\disjSet$ induce $\singleCertificate$; $\activeConstraintIndexCollection$ determines $\singleCertificate$.
			\Ensure There exists $ t \in \disjTermsIndexSet $ such that $\activeConstraintIndexSett$ is feasible for $\Qkt[\ell]$ and $\vt \neq 0$.
			\State $ \disjTermsIndexSubset[f] = \{t \in \disjTermsIndexSet : \vt \neq 0, (\mxrow{A^{\ell t}}{\activeConstraintIndexSett})^{-1} b^{\ell t}_{\activeConstraintIndexSett} \in \Qkt[\ell] \} $  \label{step:feasible_terms}  \Comment Terms meeting \Cref{lem:tighten_parameterization} Condition 1
			\State $ \singleCertificate[\barvt] \gets \singleCertificate $ \label{step:copy_certificate}
			\State $(\alpha, \beta) \gets$ \Cref{t:validate-parameterization}(\IP{\ell}, $\disjSubset, \singleCertificateSubset )$ \label{step:initial_alpha} \Comment Find cut coefficients $\alpha$
			\ForAll{$t \in \feasibleDisjTermsIndexSetk[\ell]$ such that $ A^{\MIPk t} \neq A^{\ell t} $ or $ t \notin \disjTermsIndexSubset[f] $} \Comment Terms violating \Cref{lem:tighten_parameterization}
				\State $\barvt \gets \underset{y \in \Dlta[\ell]}{\arg\max}\{y^{\T} b^{\ell t}\}$ \label{step:solve-dual} \Comment Reparameterize
			\EndFor \label{step:end_for}
			\State $(\alpha, \bar{\beta}) \gets$ \Cref{t:validate-parameterization}(\IP{\ell}, $\feasibleDisjSetk[\ell], \feasibleSingleCertificatek[\ell][\barvt])$ \label{step:tighten} \Comment Regenerate
			\State \Return $(\alpha, \bar{\beta})$
		\end{algorithmic}
	\end{algorithm}
%
	In \Cref{step:initial_alpha}, we obtain an initial valid disjunctive cut. For terms risking loss of support (\Cref{p:tight_matrix,p:tight_term,p:tight_basis}), \Cref{step:solve-dual} computes supporting certificates via \ref{Dual-lt-alpha}. Finally \Cref{step:tighten} regenerates a disjunctive-hull-supporting cut per \Cref{lem:tighten_parameterization}, as shown in \Cref{t:tighten_parameterization}.
	\begin{theorem}
		\label{t:tighten_parameterization}
		Let $k, \ell \in \MIPSet$. Let $\disjSet$ be a disjunction; $\singleCertificate$ be Farkas multipliers induced by \IP{k} and $\disjSet$; $\activeConstraintIndexCollection$ be a basis determining $\singleCertificate$; and $\{(\alpha^k, \beta^k)\}_{k \in \MIPSet}$ be the associated family of Farkas PDIs. If there exists $t \in \disjTermsIndexSet$ such that $\activeConstraintIndexSett$ is feasible for $\Qkt[\ell]$ and $\vt > 0$, then $(\alpha^\ell, \bar{\beta}^\ell)$ output from \Cref{alg:strong_param_vpcs} supports $\cl\conv(\cup_{t \in \disjTermsIndexSet} \Qkt[\ell])$.
	\end{theorem}
	Thus \textit{any} previously generated valid disjunctive inequality with a Farkas certificate can be tightened to support the disjunctive hull. 

    \section{Computational Experiments}
    \label{sec:computation_ipco}

	We evaluate implementations of \Cref{def:farkas_family} and \Cref{alg:strong_param_vpcs} on randomly perturbed MIPLIB 2017 instances, comparing against \citet{Kazachkov22vpc} and a baseline in which no disjunctive cuts are generated. All runs use standard MILP solvers on a shared server cluster; details below. Bold terms match figure legends and table headers.
	\paragraph{Test Set Generation.}
	From \textbf{Base} instances $\base$ (presolved MIPLIB-2017; solve time $\ge 10$s; $\le 5{,}000$ variables and constraints), we build \textbf{Test} sets indexed by $(k,\probdata,\theta)\in \base\times\{A,b,c\}\times\degrees$, where \textbf{Element} $=\probdata$ and \textbf{Degree} $=\theta$. Each instance $\ell$ in $(k,\probdata,\theta)$ is produced by \Cref{alg:find_perturbation}, which rotates $\probdata^\MIPk$ by angle $\theta$. For each $(k,\probdata,\theta)$ we iterate until we obtain 5 instances, hit 1{,}000 attempts, or reach 4 hours; hence test-set sizes vary by base instance. We use these test sets to approximate the diversity of possible MILP sequences while avoiding decomposition-specific or real-time–specific structure, enabling a general assessment of our methods.
	\paragraph{Methods Compared.}
	Let $d$ denote the number of \textbf{Terms} in the disjunction used for cut generation. For each $(k,\probdata,\theta)$ and instance $\ell$, we augment the root cut pool with:
	\begin{itemize}
		\item \textbf{\vpc}: $\mathcal{V}$-polyhedral disjunctive cuts \cite{Kazachkov22vpc};
		\item \textbf{\pdc}: PDIs via \Cref{def:farkas_family};
		\item \textbf{\spdc}: disjunctive-hull-supporting PDIs via \Cref{alg:strong_param_vpcs} (skipped when $\probdata=c$, since \Cref{lem:tighten_parameterization} implies \Cref{t:validate-parameterization} already supports);
		\item \textbf{\default}: default solver (no disjunctive cuts).
	\end{itemize}
	For \pdc and \spdc, we reuse the Farkas certificate $\singleCertificate$ and its determining basis $\activeConstraintIndexCollection$ and disjunction $\disjSet$ obtained when generating disjunctive cuts for \IP{k}. In what follows, we refer to the disjunction employed by \pdc and \spdc as the PDI disjunction (\PDInit) and to the disjunction employed by \vpc as the VPC disjunction (\VDInit).
	\paragraph{Software and Hardware Dependencies.}
	All experiments use Gurobi 10 with a root-node callback for adding cuts. We leverage the VPC C++ library \cite{vpc} as the implementation for \citet[Alg.~1]{Kazachkov22vpc}, which uses CBC 2.10 to extract disjunctions from partial branch-and-bound trees. Runs execute on a shared cluster (16 AMD Opteron 6128); each job uses 2 CPUs and up to 15\,GB RAM. Time limits are 3{,}600\,s for disjunctive-cut generation and 3{,}600\,s for solving.

	\subsection{Strength of Root Relaxations}
    \label{subsec:root_results}
	In this section, we show that PDIs effectively tighten root relaxations. We analyze \Cref{tab:int_gap_ipco,tab:root_gap_ipco,tab:root_time_ipco} (see \Cref{a:tables_figures}) to compare cut strength and root-node efficiency across generators. The ``Root Node'' columns of \Cref{tab:instance_counts} report the counts of unique base and test instances where root-node cut generation succeeded for all combinations of degree, term, and method.
	\begin{table}[h!]
		\centering
		\caption{Average percent integrality gap closed by disjunctive cuts alone and implied by their generating disjunctions}
		\label{tab:int_gap_ipco}
		\begin{tabular}{
			*{7}{c:}
			c
		}
			\toprule
			{Degree}
			& {Terms}
			& {Element}
			& {\VDInit}
			& {\PDInit}
			& {\vpc}
			& {\spdc}
			& {\pdc} \\
			\midrule
			\multirow{6}{*}{0.5}
			&    & matrix    & 5.96  & 2.67  & 3.38  & 1.90      & 1.90 \\
			& 4  & objective & 5.37  & 4.74  & 3.60  & \text{--} & 3.18 \\
			&    & rhs       & 8.58  & 3.69  & 3.32  & 2.23      & 2.23 \\
			\cmidrule(lr){2-8}
			&    & matrix    & 12.61 & 7.88  & 5.94  & 3.51      & 3.31 \\
			& 64 & objective & 13.70 & 12.09 & 7.88  & \text{--} & 6.21 \\
			&    & rhs       & 17.66 & 9.99  & 11.16 & 6.47      & 5.21 \\
			\cmidrule(lr){1-8}
			\multirow{6}{*}{2}
			&    & matrix    & 8.33  & 2.53  & 4.30  & 1.10      & 1.10 \\
			& 4  & objective & 6.06  & 4.35  & 4.21  & \text{--} & 2.83 \\
			&    & rhs       & 18.63 & 3.00  & 3.14  & 0.48      & 0.48 \\
			\cmidrule(lr){2-8}
			&    & matrix    & 16.26 & 5.75  & 7.43  & 1.18      & 0.44 \\
			& 64 & objective & 15.59 & 10.87 & 8.52  & \text{--} & 5.22 \\
			&    & rhs       & 27.38 & 7.31  & 9.13  & 1.16      & 1.08 \\
			\bottomrule
		\end{tabular}
	\end{table}
	\paragraph{Root LP with only disjunctive cuts (\Cref{tab:int_gap_ipco}).} As expected, each relaxation is bounded by the dual bound of its generating disjunction. \vpc gives the strongest bounds, followed by \spdc then \pdc---consistent with design: \vpc optimizes per instance; \spdc and \pdc reuse bases from prior instances, with \spdc tighter because it enforces hull support.

	Lower perturbation and larger disjunctions yield stronger parametric cuts (bases closer to optimality; tighter relaxations), but overall root-bound gains remain modest. This is expected: many disjunctive cuts are most effective deeper in the tree and are limited by the strength of the disjunctive relaxation. In our set, instances with strong relaxations were often excluded because they were too large, could not be randomly perturbed, or failed to produce non-parametric disjunctive cuts for the base instance, leaving no certificate for our methods to parameterize. We would expect the gap between \pdc and \spdc to widen with greater perturbation and larger disjunctions, but this pattern is not uniform. It does appear, along with larger root-gap gains, when the underlying disjunctive relaxation is strong \Cref{tab:example_int_gap_ipco,a:tables_figures}, suggesting the irregularities reflect limits of the relaxation rather than of the parametric methods.
%
%
	\begin{table}
		\centering
		\renewcommand{\tabcolsep}{5pt}
		\sisetup{
			table-number-alignment = center
		}
		\caption{Average percent integrality gap closed after root node processing}
		\label{tab:root_gap_ipco}
		\begin{tabular}{
			*{6}{c:}
			c
		}
			\toprule
			{Degree}
			& {Terms}
			& {Element}
			& {\vpc}
			& {\spdc}
			& {\pdc}
			& {\default} \\
			\midrule
			\multirow{6}{*}{0.5}
			&    & matrix    & 73.26 & 73.25     & 73.62 & 73.34 \\
			& 4  & objective & 61.21 & \text{--} & 61.20 & 60.07 \\
			&    & rhs       & 51.76 & 48.47     & 49.65 & 50.20 \\
			\cmidrule(lr){2-7}
			&    & matrix    & 73.58 & 73.59     & 73.40 & 73.42 \\
			& 64 & objective & 62.05 & \text{--} & 61.69 & 60.51 \\
			&    & rhs       & 53.47 & 53.02     & 52.69 & 49.72 \\
			\cmidrule(lr){1-7}
			\multirow{6}{*}{2}
			&    & matrix    & 80.24 & 80.39     & 80.41 & 80.40 \\
			& 4  & objective & 61.84 & \text{--} & 62.34 & 61.06 \\
			&    & rhs       & 52.80 & 52.45     & 52.89 & 53.77 \\
			\cmidrule(lr){2-7}
			&    & matrix    & 80.61 & 80.24     & 79.39 & 80.12 \\
			& 64 & objective & 62.68 & \text{--} & 62.39 & 61.14 \\
			&    & rhs       & 56.35 & 53.47     & 53.17 & 53.91 \\
			\bottomrule
		\end{tabular}
	\end{table}
	\paragraph{Root LP after full root processing (\Cref{tab:root_gap_ipco}).} Adding disjunctive cuts improves the root gap relative to \default, as in \citet{Kazachkov22vpc}. The strength ordering broadly persists but is less pronounced. Interactions with default root cut generation can blur differences and occasionally make nominally stronger disjunctive cuts yield a similar or weaker overall bound on this test set. As before, larger gains in root relaxation strength appear when the disjunctive relaxations themselves are stronger, as illustrated in \Cref{tab:example_root_gap_ipco} of \Cref{a:tables_figures}.
%
%
	\begin{table}[ht]
		\centering
		\renewcommand{\tabcolsep}{5pt}
		\sisetup{
			table-number-alignment = center
		}
		\caption{Average time (s) to process root node}
		\label{tab:root_time_ipco}
		\begin{tabular}{
			*{6}{c:}
			c
		}
			\toprule
			{Degree}
			& {Terms}
			& {Element}
			& {\vpc}
			& {\spdc}
			& {\pdc}
			& {\default} \\
			\midrule
			\multirow{6}{*}{0.5}
			&    & matrix    & 1.394  & 0.480     & 0.435 & 0.314 \\
			& 4  & objective & 1.893  & \text{--} & 0.601 & 0.572 \\
			&    & rhs       & 0.690  & 0.297     & 0.300 & 0.275 \\
			\cmidrule(lr){2-7}
			&    & matrix    & 11.878 & 2.135     & 1.061 & 0.336 \\
			& 64 & objective & 19.130 & \text{--} & 1.309 & 0.530 \\
			&    & rhs       & 6.428  & 0.744     & 0.549 & 0.313 \\
			\cmidrule(lr){1-7}
			\multirow{6}{*}{2}
			&    & matrix    & 2.450  & 0.377     & 0.373 & 0.296 \\
			& 4  & objective & 1.758  & \text{--} & 0.631 & 0.533 \\
			&    & rhs       & 0.794  & 0.248     & 0.234 & 0.239 \\
			\cmidrule(lr){2-7}
			&    & matrix    & 9.564  & 1.704     & 1.211 & 0.288 \\
			& 64 & objective & 19.378 & \text{--} & 1.408 & 0.540 \\
			&    & rhs       & 4.358  & 0.546     & 0.515 & 0.234 \\
			\bottomrule
		\end{tabular}
	\end{table}

    \begin{table}[htbp]
      \centering
      \caption{Unique test instance and total experiment count by setup}
      \label{tab:instance_counts}
      \begin{tabular}{l|cc|cc|cc}
        \toprule
        Element &
        \multicolumn{2}{c|}{Root Node} &
        \multicolumn{2}{c|}{Overall} &
        \multicolumn{2}{c}{Overall $>$ 20} \\
        \cmidrule(lr){2-3}\cmidrule(lr){4-5}\cmidrule(lr){6-7}
        & Base & Test & Base & Experiments & Base & Experiments \\
        \midrule
        Matrix    & 89  & 256  & 89  & 460  & 24 & 51  \\
        Objective & 106 & 578  & 106 & 1357 & 30 & 126 \\
        RHS       & 59  & 133  & 59  & 250  & 8  & 21  \\
        \bottomrule
      \end{tabular}
    \end{table}

	\paragraph{Root processing time (\Cref{tab:root_time_ipco}).}
	Runtimes reverse the strength ranking: \vpc	$>$ \spdc $>$ \pdc $>$ \default (slow to fast). This mirrors computational burden: \vpc solves cold-started LPs; \spdc solves warm-started LPs; \pdc uses fixed matrix operations; \default adds no disjunctive cuts. Processing time increases with the number of disjunctive terms, consistent with cut generation complexity scaling in disjunction size.
	\paragraph{Takeaway.} The methods form a near Pareto frontier trading root-relaxation strength for time: \vpc (strongest, slowest), followed by \spdc and \pdc, then finally \default (weakest, fastest). Because overall performance depends on both root strength and generation cost, maintaining multiple generators offers useful flexibility. The next section highlights how combining multiple cut generators with higher-term disjunctions can improve downstream solver performance.

	\subsection{Overall Performance}
    \label{subsec:bnc_results}

\begin{figure}[htbp]
		\centering
		\begin{minipage}{0.45\textwidth}
			\centering
			\includegraphics[width=\textwidth]{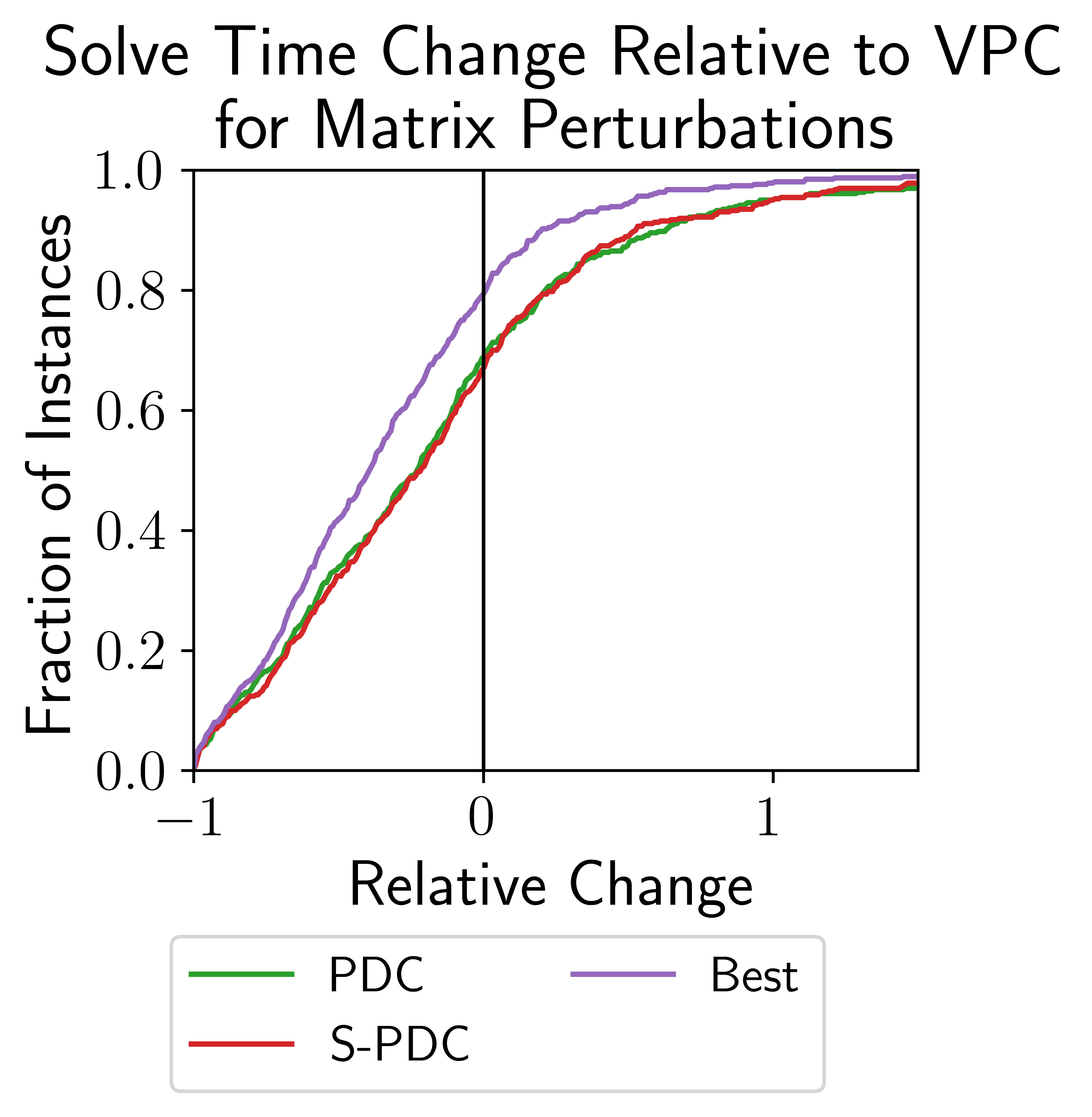}
		\end{minipage}
		\hspace{0.025\textwidth}
		\begin{minipage}{0.475\textwidth}
			\centering
			\includegraphics[width=\textwidth]{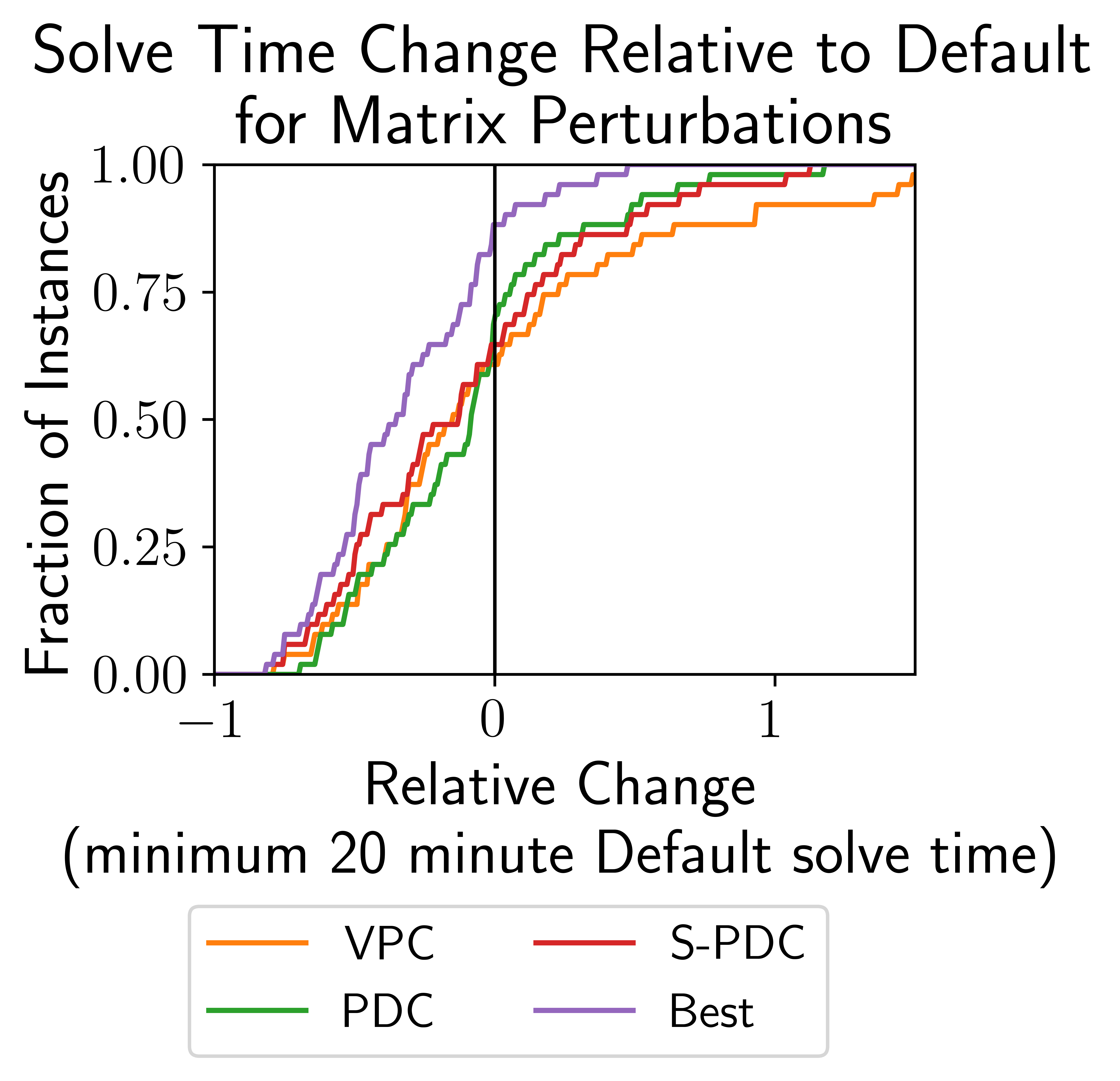}
		\end{minipage}
		\caption{For matrix perturbations, PDIs often outperform VPCs (left), and the best disjunctive cut generator almost always beats default Gurobi (right)}
		\label{fig:results_matrix}
	\end{figure}

\begin{figure}[htbp]
		\centering
		\begin{minipage}{0.45\textwidth}
			\centering
			\includegraphics[width=\textwidth]{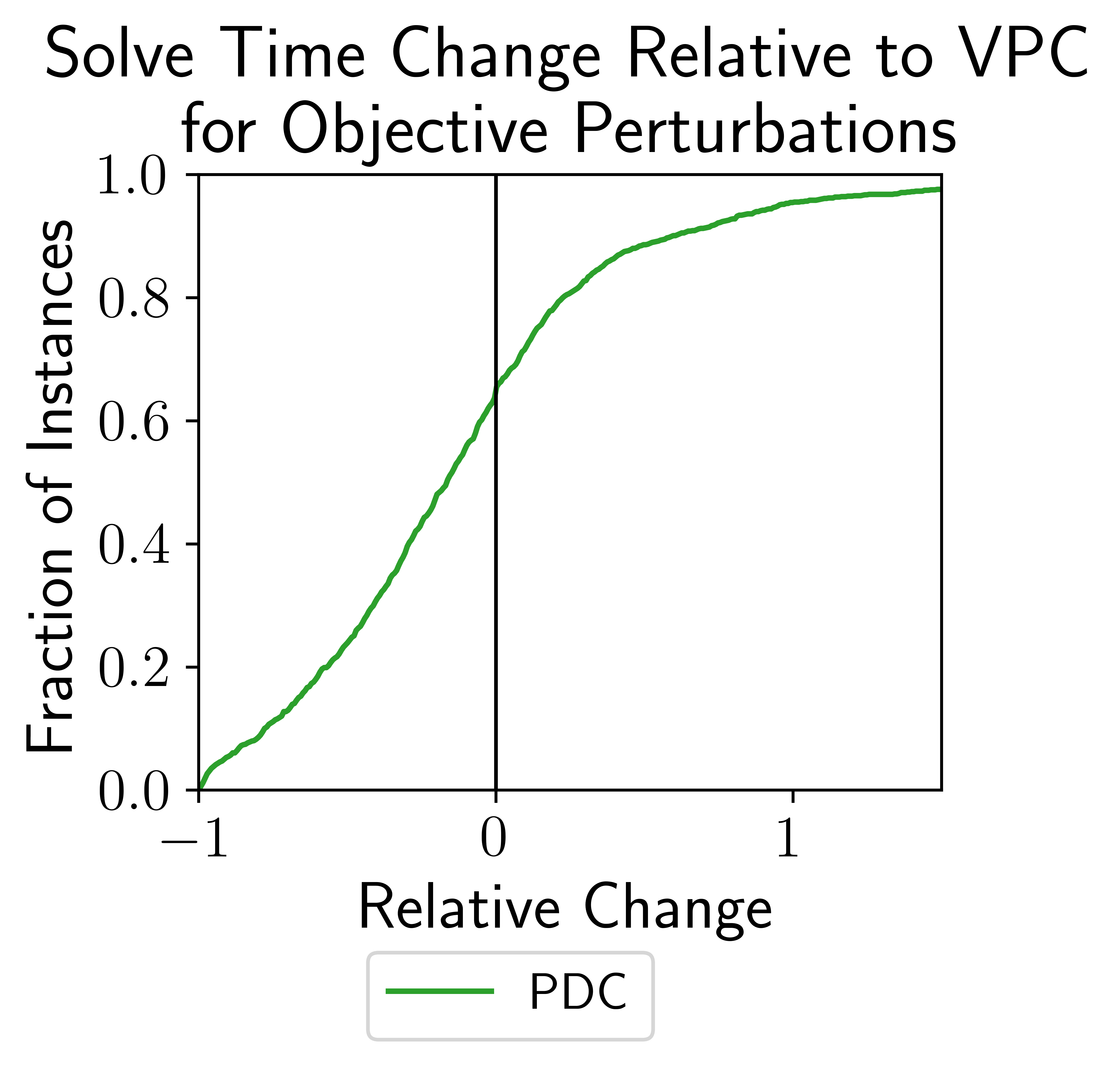}
		\end{minipage}
		\hspace{0.025\textwidth}
		\begin{minipage}{0.45\textwidth}
			\centering
			\includegraphics[width=\textwidth]{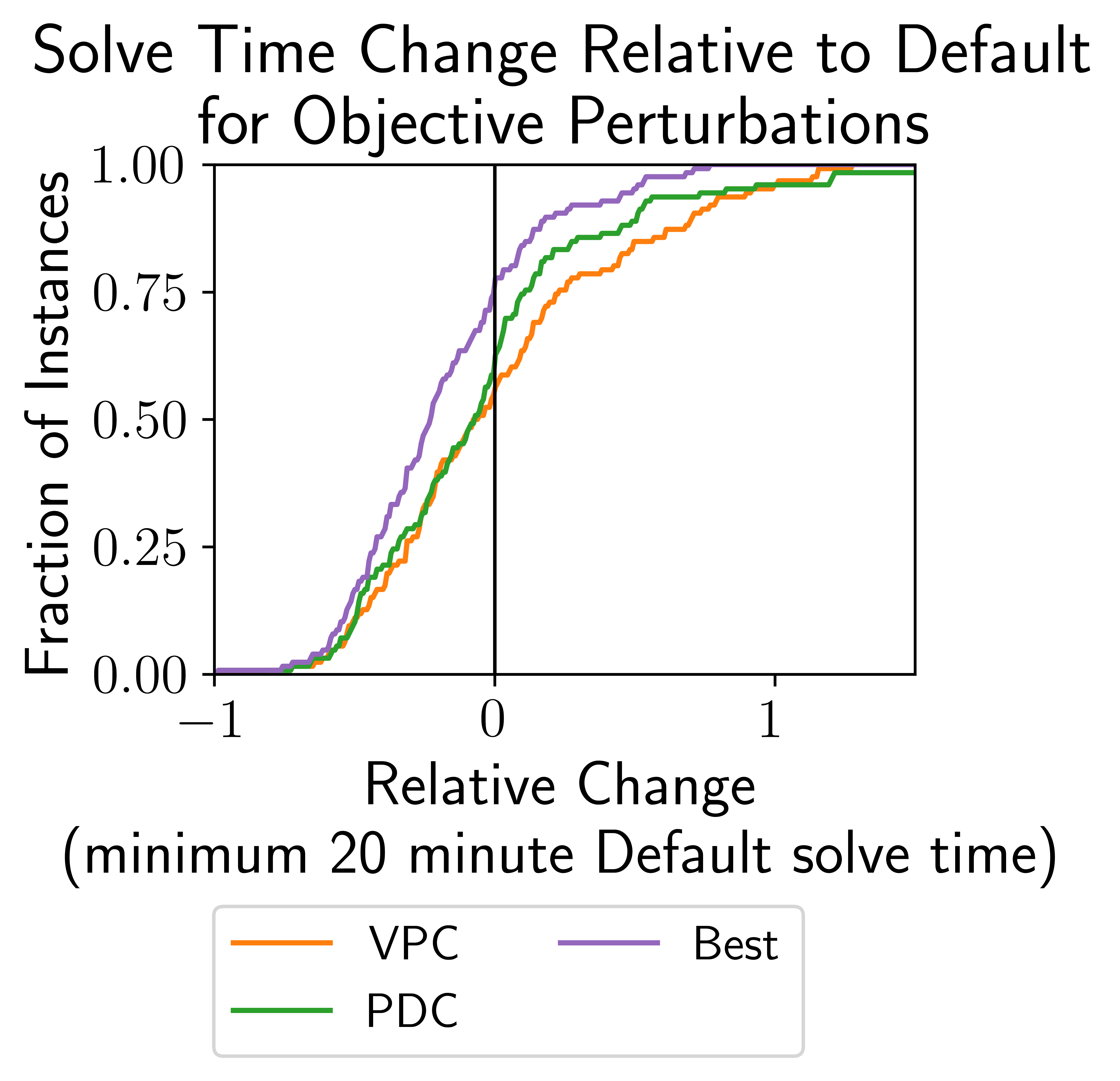}
		\end{minipage}
		\caption{PDIs often outperform VPCs (left) and default Gurobi (right) for objective perturbations}
		\label{fig:results_objective}
	\end{figure}

\begin{figure}[htbp]
		\centering
		\begin{minipage}{0.45\textwidth}
			\centering
			\includegraphics[width=\textwidth]{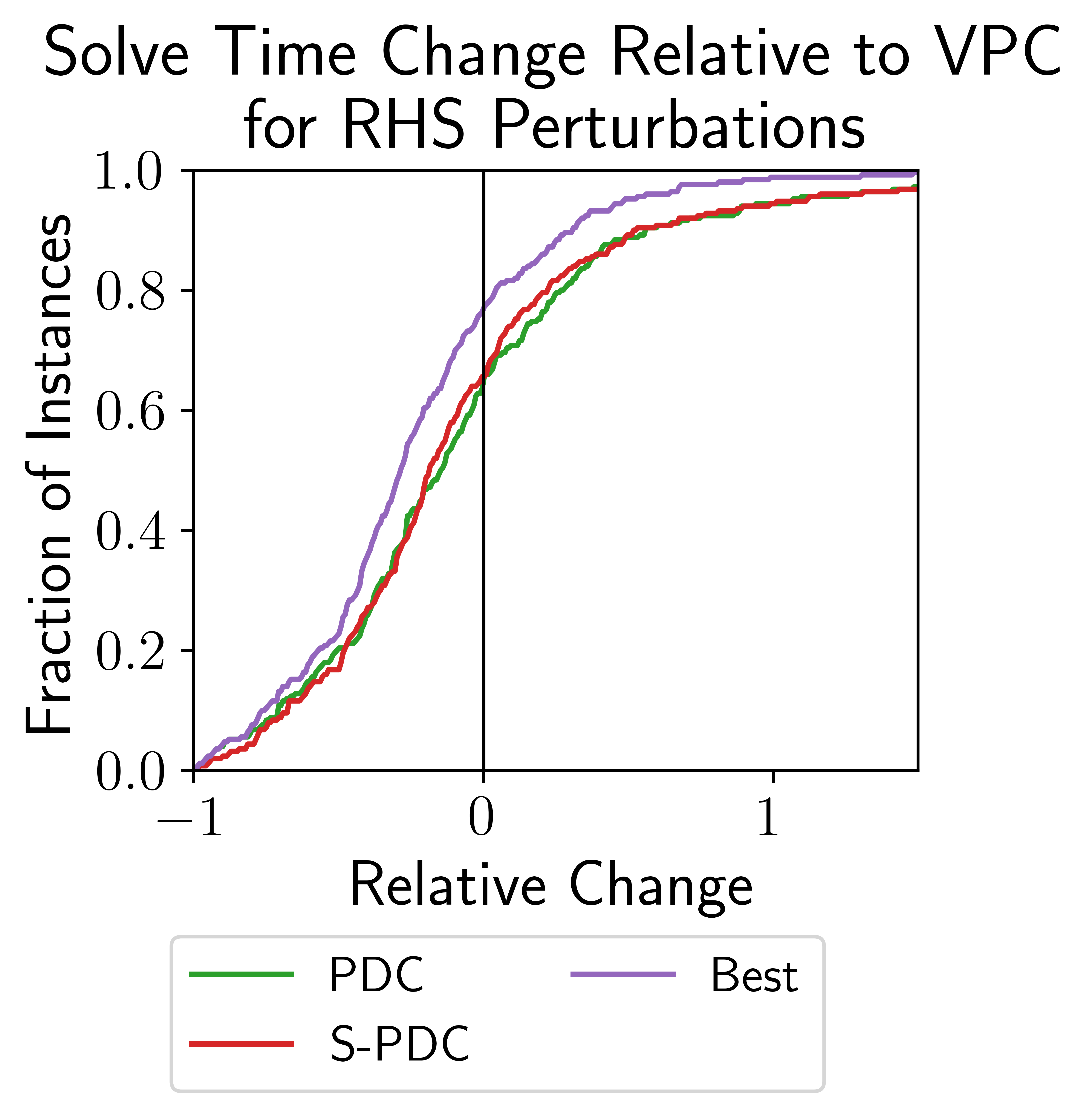}
		\end{minipage}
		\hspace{0.025\textwidth}
		\begin{minipage}{0.45\textwidth}
			\centering
			\includegraphics[width=\textwidth]{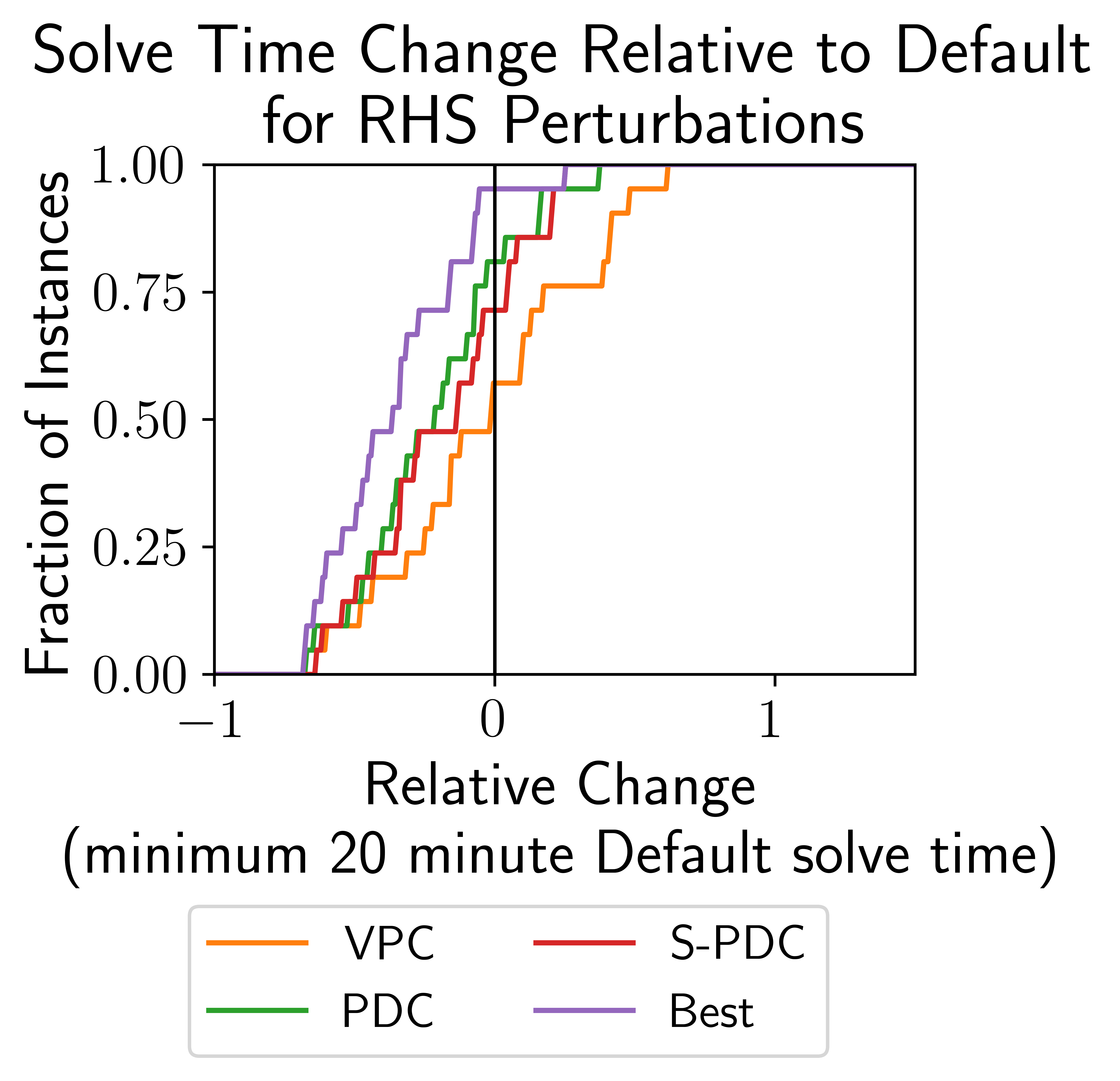}
		\end{minipage}
		\caption{PDIs often outperform VPCs and default Gurobi for right-hand side perturbations}
		\label{fig:results_rhs}
	\end{figure}

    In this section, we show that  PDIs outperform both VPCs and Gurobi with default settings.
	We analyze \Cref{fig:results_matrix,fig:results_objective,fig:results_rhs} to compare overall branch-and-cut solve time of Gurobi configured with \vpc, \spdc, \pdc, and \default.	In each figure, the left panel is a performance profile of \spdc and \pdc relative to \vpc over all test instances, where a value of -0.5 denotes a 50\% speedup and +0.5 a 50\% slowdown on instances solved by all methods. The right panel compares \vpc, \spdc, and \pdc against \default, restricted to instances where \default required at least 20 minutes. ``Best'' denotes the virtual best over all configurations. The ``Overall'' and ``Overall $>$ 20'' columns of \Cref{tab:instance_counts} report, respectively, the counts of unique base instances and total experiments that reached optimality, and the subset for which \default took $>$ 20 minutes. These need not match the count of test instances in ``Root Node'' because (a) each test instance is rerun for every degree, (b) some runs that finish root-node processing still time out before optimality, and (c) some base instances fail to yield nonparametric disjunctive cuts for certain (degree, term) pairs—so the parametric methods cannot start—whereas many runs that do receive certificates subsequently reach optimality.
	\paragraph{Matrix perturbations: relative to \vpc.}
	In \Cref{fig:results_matrix} (left), both \spdc and \pdc beat \vpc on roughly 65\% of instances. The virtual best improves on \vpc in about 80\% of cases, indicating complementary strengths between \spdc and \pdc. This aligns with \Cref{tab:int_gap_ipco,tab:root_gap_ipco,tab:root_time_ipco}: \spdc and \pdc trade a small amount of cut strength for markedly lower generation cost, yielding better end-to-end times.

	\paragraph{Matrix perturbations: relative to \default.} In \Cref{fig:results_matrix} (right), \vpc, \spdc, and \pdc each outperform \default on approximately 60–65\% of instances with solve times above 20 minutes. The virtual best wins on nearly 90\% of these cases, with about half solving at least twice as fast. This points to a usual branch-and-bound tradeoff: stronger root cuts increase root processing time, but they tighten the relaxation globally. On harder, long-running instances, implicit pruning from the tighter relaxation can accumulate and the upfront cost is amortized, creating opportunity to lower total solve time.

	\paragraph{Objective and right-hand side perturbations.} The same patterns hold for objective-coefficient and right-hand-side perturbations: \spdc and \pdc compare favorably to \vpc on the left panels, and all three disjunctive generators improve on \default for long-running instances on the right panels (\Cref{fig:results_objective,fig:results_rhs}, \Cref{a:tables_figures}). These results reinforce that the observed trends are robust across perturbation types and not specific to the constraint matrix.

    \section{Conclusion}
    \label{sec:conclusion_ipco}

	This paper provides theory on how to generate PDIs that remain valid (and supporting) across related MILPs without solving a cut-generating LP; dual solutions of the tightening LPs certify support for individual terms and, when fed into \Cref{t:validate-parameterization}, yield cuts that provably support the disjunctive hull. In our experiments, strengthened and unstrengthened PDIs tighten the default root relaxation while retaining much of the speed advantage over \citet{Kazachkov22vpc}. Within branch and cut, our best parametric generator consistently outperforms \citet{Kazachkov22vpc} in total solve time, and the strongest disjunctive generator beats default Gurobi on most instances exceeding 20 minutes.
    
    Open questions remain, notably how to choose the most effective inequality per instance; under matrix perturbations, deriving the supporting certificate from the original (nonparametric) cut may further help.
    
    Extensions include learning-guided selection of generators, comparing our disjunction-based warm start with an explicit MILP warm start using the same disjunction, and applying our strengthening routine to other costly families (e.g., bilevel interdiction cuts) via \Cref{alg:strong_param_vpcs}.
    
    Methodologically, we provide a pathway to parameterize disjunctive inequalities and recover hull support as needed. Practically, we offer guidance for warm-starting the dual side---especially the cut pool---offering reduced solve times in sequential and decomposed settings and highlighting gains from exploiting structure across sequences of related problems.

    \newpage

    \bibliographystyle{splncs04nat}
    \bibliography{references/reference}

    \newpage

    \appendix

	\section{Proofs}
	\label{a:proofs}

	\subsubsection{Proof of \Cref{t:validate-parameterization}}
	\label{proof:validate_parameterization}

	\begin{proof}
        Let
        $
        \gamma^t \defeq
        v^t A^{kt}
        $
        and
        $
        \gamma_0^t \defeq v^t b^{kt}
        $
        for all $t \in \disjTermsIndexSet$.
        We will show that $\alpha^\T x \ge \beta$ is valid for $\Qkt$ for a fixed (arbitrary) index $t \in \disjTermsIndexSet$.
        We have that
        $\gamma^t x \ge \gamma^t_0$
        is valid for $\Qkt$.
        Then, if $x \in \Qkt$, since we have nonnegativity on the variables ($\Qkt \subseteq \nonnegreals^n$),
        it follows that,
        \begin{equation*}
            \alpha^\T x
            =
            \sum_{j \in [n]} \max_{t' \in \disjTermsIndexSet} \{\gamma^{t'}_j\} x_j
            \ge
            \sum_{j \in [n]} \gamma^t_j x_j
            \ge
            \gamma^t_0
            \ge
            \min_{t' \in \disjTermsIndexSet} \{\gamma^{t'}_0\}
            = \beta.
        \end{equation*}
        Hence, $\alpha^\T x \ge \beta$ is valid for $\bigcup_{t\in\disjTermsIndexSet}\Qkt \supseteq \PIk \cap \PDk$.
    \end{proof}

	\subsubsection{Proof of \Cref{lem:tighten_parameterization}}
	\label{proof:tighten_parameterization_lemma}

	\begin{proof}
        Let $t' = \underset{t \in \disjTermsIndexSet}{\arg\min}\{v^t b^{\ell t}\}$. Our hypothesis provides $v^t A^{\ell t} = v^t A^{k t} = \alpha$ for all $t \in \disjTermsIndexSet$, which implies that $\underset{t \in \disjTermsIndexSet}{\max}\{v^t A^{\ell t}\} = v^{t'} A^{\ell t'}$. Applying \Cref{t:validate-parameterization} yields $(\alpha, \beta) = (\underset{t \in \disjTermsIndexSet}{\max}\{v^t A^{\ell t}\}, \underset{t \in \disjTermsIndexSet}{\min}\{v^t b^{\ell t}\}) = (v^{t'} A^{\ell t'}, v^{t'} b^{\ell t'})$, a valid inequality for $\cl\conv(\cup_{t \in \disjTermsIndexSet} \Qkt[\ell])$. Assumptions on $\activeConstraintIndexSett$ imply that $ \alpha = v^{t'}_{\activeConstraintIndexSett} \mxrow{A^{\ell t'}}{\activeConstraintIndexSett} $, $ v^{t'}_{\activeConstraintIndexSett} = \alpha (\mxrow{A^{\ell t'}}{\activeConstraintIndexSett})^{-1} $, $ v^{t'}_{\activeConstraintIndexSett} b^{\ell t'}_{\activeConstraintIndexSett} = \alpha (\mxrow{A^{\ell t'}}{\activeConstraintIndexSett})^{-1} b^{\ell t'}_{\activeConstraintIndexSett} $, and $ (\mxrow{A^{\ell t'}}{\activeConstraintIndexSett})^{-1} b^{\ell t'}_{\activeConstraintIndexSett} \in \Qkt[\ell] $. Since $ \beta = v^{t'}_{\activeConstraintIndexSett} b^{\ell t'}_{\activeConstraintIndexSett} $, $(\alpha, \beta)$ supports $\Qktprime[\ell]$. Given $\PDk[\ell] = \cl\conv(\cup_{t \in \disjTermsIndexSet} \Qkt[\ell])$ and is convex, $(\alpha, \beta)$ supports it, too.
    \end{proof}

	\subsubsection{Proof of \Cref{lem:supporting_certificate_exists}}
	\label{proof:supporting_certificate_exists}

	\begin{proof}
		Let $\opt{x} \in \arg\min_{x \in \Qkt[\ell]}\{\alpha^{\T} x\}$ and $\barvt \in \arg\max_{\vt \in \Dlta}\{\vt b^{\ell t}\}$. Since $\opt{x}$ optimal for \ref{LP-lt-alpha}, we have that $\alpha^{\T} x \geq \alpha^{\T}\opt{x}$ for all $x \in \Qkt[\ell]$, which means that $\alpha, \alpha^{\T}\opt{x}$ supports $\Qkt$. Since $\barvt$ is feasible for \ref{Dual-lt-alpha}, we have that $\barvt A^{\ell t} = \alpha$. By stong duality, $\barvt b^{\ell t} = \alpha^{\T}\opt{x}$, which means that $\barvt$ is a certificate for the cut $(\alpha, \alpha^{\T}\opt{x})$ and is thus supporting.
	\end{proof}

	\subsubsection{Proof of \Cref{cor:supporting_certificate_unique}}
	\label{proof:supporting_certificate_unique}

	\begin{proof}
		\Cref{lem:supporting_certificate_exists} provides existence of $\opt{x}$ and $\barvt$. Since there does not exist $ i \in [\numRowsQkt] $ such that $ \alpha = \mxrow{A^{\MIPk t}}{i} $, pivoting from $\opt{x}$ to another $x \in \Qkt$ would violate dual feasibility. Therefore, $\opt{x} \in \Qkt$ and $\opt{y} \in \Dlta$ such that $ \alpha^{\T} \opt{x} = \opt{y}^{\T} b^{\MIPk t} $ are unique.
	\end{proof}

	\subsubsection{Proof of \Cref{t:tighten_parameterization}}
	\label{proof:tighten_parameterization_theorem}

	\begin{proof}
		Define $\disjTermsIndexSubset[f]$ as in \Cref{step:feasible_terms}, which is nonempty due to the existence of $t \in \disjTermsIndexSet$ such that $\activeConstraintIndexSett$ is feasible for $\Qkt[\ell]$ and $\vt > 0$. As prescribed in \Cref{step:copy_certificate}, let $\singleCertificate[\barvt] = \singleCertificate$, and, additionally, let $\activeConstraintIndexCollection[\barA] = \activeConstraintIndexCollection$. By \Cref{lem:tighten_parameterization}, \Cref{step:initial_alpha} provides $\alpha = \barvt A^{\ell t}$ for all $t \in \disjTermsIndexSubset[c]$. Consider $t \in \feasibleDisjTermsIndexSetk[\ell] \setminus \disjTermsIndexSubset[c]$. Update $\barAt$ such that $\mxrow{A^{\ell t}}{\barAt} \bar{x}^t = b^{\ell t}_{\barAt}$ for $\bar{x}^t \in \arg\max_{x \in \Qkt[\ell]}\{\alpha^{\T} x\}$. \Cref{lem:supporting_certificate_exists} proves in \Cref{step:solve-dual} that  $\barvt \in \Rplus^{1 \times n}$ and $\alpha = \barvt A^{\ell t}$. Since $\bar{x}^t$ and $\barvt$ are dual to each other, we have that $\{i \in [\numRowsQkt] : \barvt_i > 0\} \subseteq \barAt$. Therefore, we have that \IP{\ell} and $\feasibleDisjSetk[\ell]$ induce $\feasibleSingleCertificatek[\ell][\barvt]$, $\AtCollection[\barA]$ determines $\feasibleSingleCertificatek[\ell][\barvt]$ and, for all $t \in \feasibleDisjTermsIndexSetk[\ell]$, $\barAt$ is feasible for $\Qkt[\ell]$. thus, \Cref{lem:tighten_parameterization} provides $(\alpha, \bar{\beta})$ in \Cref{step:tighten} supports $\cl\conv(\cup_{t \in \feasibleDisjTermsIndexSetk[\ell]} \Qkt[\ell]) = \cl\conv(\cup_{t \in \disjTermsIndexSet} \Qkt[\ell])$.
	\end{proof}

	\section{Algorithms}

	\begin{algorithm}
        \caption{Find Degree}
        \label{alg:find_degree}
        \begin{algorithmic}[1]
            \Require $\probdata^\MIPk, \vnew$
            \State angle\_difference $\gets \cos^{-1} \left(\frac{\probdata^\MIPk \cdot \vnew}{\|\probdata^\MIPk\| \|\vnew\|}\right)$  \Comment{Angle between vectors}
            \State norm\_difference $\gets \left\lvert \frac{\Vert \probdata^\MIPk \Vert - \Vert \vnew \Vert}{\Vert \probdata^\MIPk \Vert} \right\rvert$  \Comment{Relative change in vector norms}
            \State \Return max$\{\text{angle\_difference, norm\_difference}\}$
        \end{algorithmic}
    \end{algorithm}

    \begin{algorithm}
        \caption{Find Perturbation}
        \label{alg:find_perturbation}
        \begin{algorithmic}[1]
            \Require $\probdata^\MIPk, \theta$ \Comment{starting vector, desired perturbation}
            \State $\probdata^\MIPk \gets \toVector(\probdata^\MIPk)$ if $\probdata^\MIPk \in \R^{\numRowsP \times n}$ else $\probdata^\MIPk$   \Comment{Flatten to vector if needed}
            \State $\vfinal, \epsilon \gets \text{NULL}, 1$
            \While{$\epsilon \geq 10^{-6}$ and $\vfinal ==$  NULL}
                \Comment{until tolerance met or perturbed vector found}
                \State $\vnew, \vprev \gets \probdata^\MIPk, \text{NULL}$
                \While{\nameref{alg:find_degree}$(\probdata^\MIPk, \vnew) < \theta$}
                    \Comment{Until desired perturbation reached}
                    \State $\vprev \gets \vnew$
                    \State $i \gets \random([|u|])$ \Comment{Select random element}
                    \State $\vnew[i] \gets \vnew[i] + \random([-\epsilon, \epsilon])$  \Comment{Randomly perturb by $\pm \epsilon$}
                \EndWhile
                \If{$\vprev \neq \probdata^\MIPk$}
                    \Comment{If an iteration produced a vector within the desired perturbation}
                    \State $\vfinal \gets \vprev$
                \EndIf
                \State $\epsilon \gets \epsilon / 2$
            \EndWhile
            \State \Return $\toMatrix(\vfinal)$ if $\vfinal \in \R^{\numRowsP n}$ else $\vfinal$  \Comment{Return as matrix if needed}
        \end{algorithmic}
    \end{algorithm}

	\section{Additional Tables and Figures}
	\label{a:tables_figures}

\begin{figure}[htbp]
		\centering
		\includegraphics[width=\textwidth]{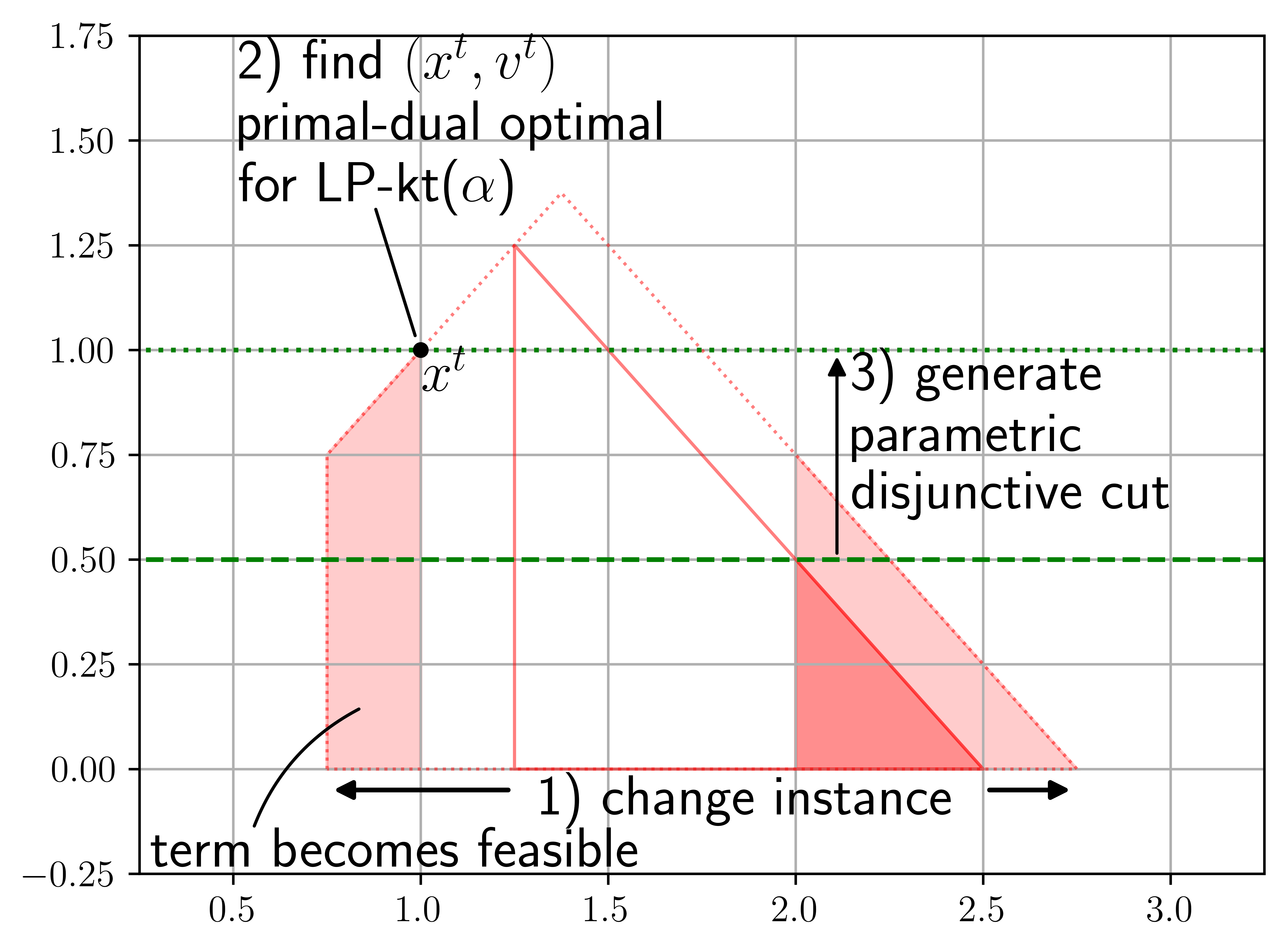}
		\caption{For perturbation-induced feasible terms}
		\label{p:tight_term}
	\end{figure}

	\Cref{p:tight_term} illustrates the second failure mode of \Cref{lem:tighten_parameterization} Condition 1: the PDI may not support the disjunctive hull when perturbation induces feasibility for a disjunctive term. As mentioned in \Cref{sec:background_ipco}, $\vt = 0$ for initially infeasible disjunctive terms, necessitating calculation of a supporting certificate.

	\Cref{tab:example_int_gap_ipco,tab:example_root_gap_ipco,tab:example_root_time_ipco} highlight the strength differences of our disjunctive cut generators compared to their generating disjunctions and solver defaults, which are made possible by the presence of a strong disjunctive relaxation.

	\begin{table}
		\centering
		\renewcommand{\tabcolsep}{5pt}
		\sisetup{
			table-number-alignment = center
		}
		\caption{Average percent integrality gap closed by disjunctive cuts alone and implied by their generating disjunctions, grouped by degree of perturbation, size of disjunction, and element perturbed for test sets generated from \texttt{bm23.mps}}
		\label{tab:example_int_gap_ipco}
		\begin{tabular}{
			*{7}{c:}
			c
		}
			\toprule
			{Degree}
			& {Terms}
			& {Element}
			& {\VDInit}
			& {\PDInit}
			& {\vpc}
			& {\spdc}
			& {\pdc} \\
			\midrule
			\multirow{6}{*}{0.5}
			&    & matrix    & 14.78 & 14.78 & 14.78 & 14.63     & 14.62 \\
			& 4  & objective & 14.63 & 14.63 & 14.63 & \text{--} & 14.63 \\
			&    & rhs       & 13.85 & 14.60 & 13.85 & 14.22     & 14.22 \\
			\cmidrule(lr){2-8}
			&    & matrix    & 71.00 & 71.35 & 69.05 & 67.23     & 66.46 \\
			& 64 & objective & 68.83 & 71.48 & 66.74 & \text{--} & 69.37 \\
			&    & rhs       & 69.03 & 70.93 & 66.91 & 68.17     & 64.77 \\
			\cmidrule(lr){1-8}
			\multirow{6}{*}{2}
			&    & matrix    & 13.38 & 14.33 & 13.22 & 13.00     & 12.90 \\
			& 4  & objective & 14.79 & 14.79 & 14.74 & \text{--} & 14.76 \\
			&    & rhs       & 11.77 & 13.27 & 11.61 & 12.16     & 12.16 \\
			\cmidrule(lr){2-8}
			&    & matrix    & 70.44 & 71.38 & 67.68 & 49.42     & 28.17 \\
			& 64 & objective & 66.90 & 70.41 & 64.82 & \text{--} & 67.60 \\
			&    & rhs       & 69.49 & 63.82 & 66.88 & 60.35     & 44.03 \\
			\bottomrule
		\end{tabular}
	\end{table}

	\Cref{tab:example_int_gap_ipco} mirrors \Cref{tab:int_gap_ipco}, reporting the percentage of integrality gap closed over 20 random perturbations of \texttt{bm23.mps}. In this case the separation of \spdc from \pdc is clearer for larger disjunctions and higher perturbation levels---precisely where parameterization has more scope to weaken. We also see the expected monotone effects: more terms and smaller perturbations yield stronger parametric cuts, and disjunctive cuts are strong overall. Apparent exceptions arise when the underlying disjunctive relaxation is stronger for the PDI disjunction than for the VPC disjunction (e.g., \PDInit $>$ \VDInit for most right-hand-side perturbations), in which case \spdc and \pdc can surpass \vpc.

	\begin{table}
		\centering
		\renewcommand{\tabcolsep}{5pt}
		\sisetup{
			table-number-alignment = center
		}
		\caption{Average percent integrality gap closed after root node processing, grouped by degree of perturbation, size of disjunction, and element perturbed for test sets generated from \texttt{bm23.mps}}
		\label{tab:example_root_gap_ipco}
		\begin{tabular}{
			*{6}{c:}
			c
		}
			\toprule
			{Degree}
			& {Terms}
			& {Element}
			& {\vpc}
			& {\spdc}
			& {\pdc}
			& {\default} \\
			\midrule
			\multirow{6}{*}{0.5}
			&    & matrix    & 42.63 & 40.62     & 40.48 & 42.16 \\
			& 4  & objective & 39.18 & \text{--} & 40.89 & 39.23 \\
			&    & rhs       & 43.19 & 41.49     & 42.00 & 43.87 \\
			\cmidrule(lr){2-7}
			&    & matrix    & 70.30 & 69.06     & 68.57 & 42.13 \\
			& 64 & objective & 68.40 & \text{--} & 70.52 & 39.38 \\
			&    & rhs       & 68.69 & 69.44     & 68.08 & 42.63 \\
			\cmidrule(lr){1-7}
			\multirow{6}{*}{2}
			&    & matrix    & 39.92 & 40.10     & 41.11 & 38.71 \\
			& 4  & objective & 40.16 & \text{--} & 39.54 & 40.02 \\
			&    & rhs       & 42.92 & 42.49     & 41.60 & 41.16 \\
			\cmidrule(lr){2-7}
			&    & matrix    & 69.56 & 56.42     & 48.44 & 38.82 \\
			& 64 & objective & 67.09 & \text{--} & 69.28 & 40.01 \\
			&    & rhs       & 69.34 & 62.25     & 56.53 & 41.16 \\
			\bottomrule
		\end{tabular}
	\end{table}

	\Cref{tab:example_root_gap_ipco} parallels \Cref{tab:root_gap_ipco}, reporting the percentage of integrality gap closed by the LP relaxation after root processing for strong disjunctive relaxations generated from random perturbations of \texttt{bm23.mps}. Given the strength of the relaxation and the clear separation among generators in \Cref{tab:example_int_gap_ipco}, the differences between \vpc, \spdc, \pdc, and \default are more pronounced here---offering a cleaner extension of the patterns seen at the ``cuts only'' level than in the aggregate tables. Deviations from monotonicity arise when the underlying disjunctive relaxation favors the PDI disjunction over the VPC disjunction (e.g., \PDInit $>$ \VDInit) and for 4-term disjunctions where disjunctive cuts fail to tighten the relaxation near the optimal region under the default cuts. Overall, the expected trends persist: more terms and smaller perturbations yield larger root-gap reductions, and stronger generators deliver greater improvements.

	\begin{table}[ht]
		\centering
		\renewcommand{\tabcolsep}{5pt}
		\sisetup{
			table-number-alignment = center
		}
		\caption{Time to process root node, grouped by degree of perturbation, size of disjunction, and element perturbed for test sets generated from \texttt{bm23.mps}}
		\label{tab:example_root_time_ipco}
		\begin{tabular}{
			*{6}{c:}
			c
		}
			\toprule
			{Degree}
			& {Terms}
			& {Element}
			& {\vpc}
			& {\spdc}
			& {\pdc}
			& {\default} \\
			\midrule
			\multirow{6}{*}{0.5}
			&    & matrix    & 0.227 & 0.103     & 0.113 & 0.134 \\
			& 4  & objective & 0.200 & \text{--} & 0.122 & 0.105 \\
			&    & rhs       & 0.213 & 0.122     & 0.124 & 0.162 \\
			\cmidrule(lr){2-7}
			&    & matrix    & 1.629 & 0.405     & 0.176 & 0.131 \\
			& 64 & objective & 1.661 & \text{--} & 0.173 & 0.103 \\
			&    & rhs       & 1.629 & 0.214     & 0.174 & 0.127 \\
			\cmidrule(lr){1-7}
			\multirow{6}{*}{2}
			&    & matrix    & 0.219 & 0.143     & 0.126 & 0.102 \\
			& 4  & objective & 0.209 & \text{--} & 0.119 & 0.124 \\
			&    & rhs       & 0.227 & 0.135     & 0.118 & 0.105 \\
			\cmidrule(lr){2-7}
			&    & matrix    & 1.702 & 0.548     & 0.214 & 0.104 \\
			& 64 & objective & 1.633 & \text{--} & 0.176 & 0.118 \\
			&    & rhs       & 1.589 & 0.364     & 0.198 & 0.106 \\
			\bottomrule
		\end{tabular}
	\end{table}

	\Cref{tab:example_root_time_ipco} parallels \Cref{tab:root_time_ipco}, reporting root processing times under strong disjunctive relaxations from random perturbations of \texttt{bm23.mps}. Since cut-generation complexity does not depend on the strength of the disjunctive relaxation, the ordering of runtimes is essentially unchanged: \vpc is slowest, followed by \spdc, then \pdc, with \default fastest. These results reinforce the expected trade-off between generator strength and processing time, with a clear negative association between strength and speed.

\end{document}